\documentclass[12pt]{article}
\title{Root games on Grassmannians}
\author{Kevin Purbhoo 
\footnote{Research partially supported by an NSERC scholarship. }\\
University of British Columbia\\
{\normalsize \tt kevinp@math.ubc.ca}}

\usepackage[mathscr]{eucal}
\usepackage{latexsym}
\usepackage{amssymb}
\usepackage{amsthm}
\usepackage{amsmath}
\usepackage{epsfig}

\newcommand{\CC}{\mathbb{C}}

\newcommand{\RR}{\mathbb{R}}

\newcommand{\skewtabl}{\lambda_{\bar 3}/\lambda_2}
\newcommand{\squares}{\mathcal{S}}
\newcommand{\Udest}{U_{\overline{23}}}
\newcommand{\ra}{\mathrm{a}}
\newcommand{\rb}{\mathrm{b}}
\newcommand{\ten}{\text{\small $1\!0$}}
\newcommand{\scriptten}{\text{$1\!0$}}

\newtheorem{theorem}{Theorem}
\newtheorem{lemma}{Lemma}[section]
\newtheorem{example}{Example}[section]
\newtheorem{definition}[example]{Definition}
\newtheorem{remark}[example]{Remark}

\newtheorem{proposition}[lemma]{Proposition}
\newtheorem{corollary}[lemma]{Corollary}

\begin{document}
\maketitle

\begin{abstract}
We recall the {\em root game}, introduced in \cite{P}, which
gives a fairly powerful sufficient condition for non-vanishing of Schubert 
calculus on a generalised flag manifold $G/B$.  
We show that it gives a {\em necessary and sufficient} rule for non-vanishing 
of Schubert calculus on Grassmannians.  In particular, a Littlewood-Richardson
number is non-zero if and only if it is possible 
to win the corresponding root game.  More generally, the rule can be used 
to determine whether or not a product of several Schubert classes on 
$Gr_l(\CC^n)$ is non-zero in a manifestly symmetric way.  Finally, we give 
a geometric interpretation of root games for Grassmannian Schubert 
problems.
\end{abstract}

%
%%%%%%%%%%%%%%%%%%%%%%%%%%%%%%%%%%%%%%%%%%%%%%%%%%%%%%%%%%%%%%
%

\section{Prior work}
\label{sec:prior}

In \cite{P} we introduced the {\em root game}, a combinatorial game which 
can often determine whether or not a given Schubert structure constant is zero
in the cohomology ring of a generalised flag manifold $G/B$.  
Our goal in this paper is to
strengthen our earlier results in the case
where the group $G$ is $GL(n, \CC)$ and the Schubert classes are pulled back 
from a Grassmannian.
We begin by recalling the root game for Schubert intersection numbers
on the (ordinary) flag manifold $Fl(n)$.  Apart from the root game,
most of the relevant background material for this paper can be found
in \cite{F1, FG}.

Let $G=GL(n)$. Let $B$ and $B_-$ denote the Borel 
subgroups of upper and lower 
triangular matrices respectively, and $T = B \cap B_-$ the standard
maximal torus, consisting of invertible diagonal matrices.

Recall that for each element of the symmetric group $\pi \in S_n$, there 
is a corresponding $T$-fixed point $\pi B$ on the flag manifold 
$Fl(n) = G/B$ (here we view $\pi$ as an element of $GL(n)$ via 
the standard representation of $S_n$), 
and an associated Schubert variety $X_\pi = \overline{B_- \cdot \pi B}$.  
We denote its cohomology class in $H^*(Fl(n))$ by $[X_\pi]$.  

Our convention will be to write all permutations in one line notation
$$ \pi = \pi(1) \pi(2) \ldots \pi(n).$$
If $1 \in S_n$ denotes the identity element, and 
$w_0 = n \ldots 3 2 1 \in S_n$ is the long word, then the 
Schubert class $[X_1]$ 
is the identity element in $H^*(G/B)$, and 
$[X_{w_0}] \in H^\text{top}(G/B)$ is the class of a point.  
In general $[X_\pi]$ is a class of degree $2\ell(\pi)$, where $\ell(\pi)$
denotes the length of $\pi$.

For $\pi_1, \ldots, \pi_m \in S_n$, 
the Schubert intersection number
\begin{equation}
\label{eqn:intnumber}
\int_{Fl(n)} [X_{\pi_1}] \cdots [X_{\pi_m}]
\end{equation}
is always a non-negative integer.  The root game attempts to determine
whether this number is strictly positive.

The game is played on a set of squares $\squares = 
\{S_{ij}\ |\ 1\leq i < j \leq n\}$.
In our diagrams, we will
arrange the squares $S_{ij}$ in an array, 
where $i$ is the row index and $j$ is the column index.
In each square we allow tokens to appear.  
Each token has a label $k \in \{1, \ldots, m\}$, and no square may ever contain 
two tokens with the same label.
A token labelled $k$ is called a $k$-token, and we write $k \in S_{ij}$
if a $k$-token appears in square $S_{ij}$.

A position in the game is specified by two pieces of data:
\begin{itemize}
\item The configuration of the tokens. Formally this is a map $\tau$
from $\squares$ to subsets of $\{1, \ldots, m\}$, and our notation
$k \in S$ is shorthand for $k \in \tau(S)$; however, in this
paper we will wish to think of each token as a physical 
object which can be moved from square to square.
\item A partition of the set of squares 
$\squares = R_1 \sqcup \cdots \sqcup R_r$.
Each $R_i$ is called a {\bf region}.
\end{itemize}

The {\bf initial position} of the game is as follows: there is a single region
$R_1 = \squares$, and for 
$i<j$, a $k$-token appears in square $S_{ij}$ 
if and only if $\pi_k(i) > \pi_k(j)$.

From the initial position we move the tokens in the manner prescribed
in the next paragraph.
However, before each move we have the option of {\em splitting} regions 
into multiple regions.  We define an {\bf ideal subset} of the squares 
to be a set $A \subset \squares$ with the property that 
if $S_{ij} \in A$, $i' \geq i$ and $j' \geq j$ then 
$S_{i'j'} \in A$.  Given an ideal subset
of the squares we can break up a region $R$ into two regions:
$R \cap A$ and $R \setminus A$.  We call this {\bf splitting $R$
along $A$}, and we may repeat the process as many times as desired.

\begin{figure}[htbp]
  \begin{center}
    \epsfig{file=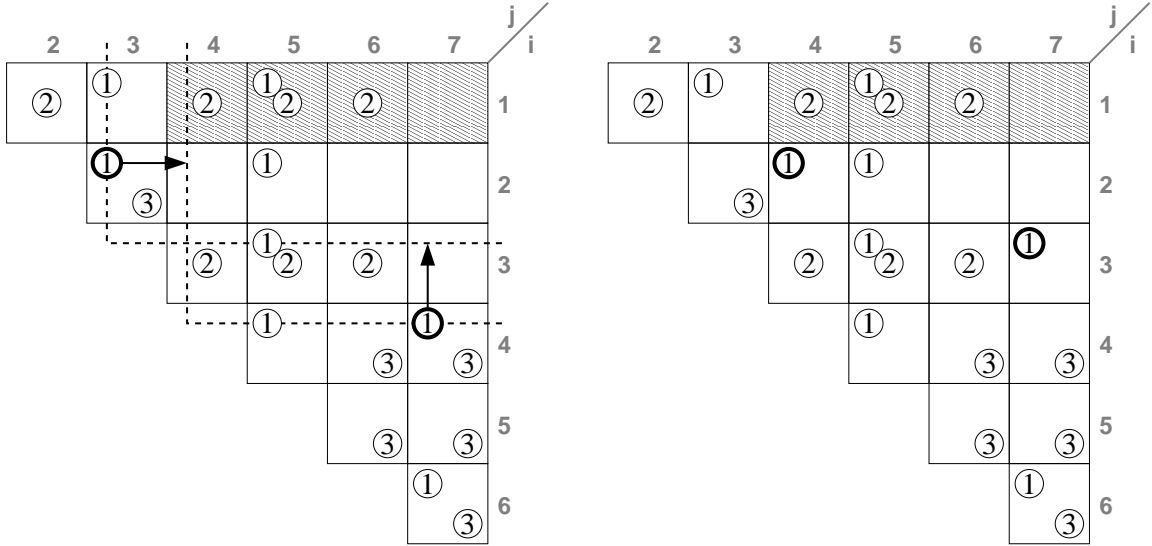,width=6in}
    \caption{On the left is the initial position of the root
 game for $\pi_1 = 3426175$, $\pi_2 = 5162347$, $\pi_3 = 1326754$,
 after splitting into two regions.  On the right is the position
 of the tokens after the move for $R =  \text{the unshaded region}$,
 $k = 1$, and $(i,j) = (3,4)$.  Note that the rows are indexed
 by $i=1, \ldots, 6$, and the columns by $j=2,\ldots, 7$, 
 since $1 \leq i < j \leq 7$.}
    \label{fig:flaggame}
  \end{center}
\end{figure}

A {\bf move} is specified by a region $R$, a token label 
$k \in \{1,\ldots,m\}$ and a pair $(i,j)$ with $1 \leq i < j \leq n$.
After choosing these data, we move tokens as follows:
\begin{itemize}
\item
For every $h$ with $j<h \leq n$, if $S_{jh}$ and $S_{ih}$ are both in $R$
and a $k$-token appears 
in $S_{jh}$ but not in $S_{ih}$, we move the $k$-token 
from $S_{jh}$ to $S_{ih}$;
\item 
For every $h$ with $1 \leq h < i$, if $S_{hi}$ and $S_{hj}$ 
are both in $R$ and a $k$-token appears in 
$S_{hi}$ but not in $S_{hj}$, we move the $k$-token 
from $S_{hi}$ to $S_{hj}$.
\end{itemize}
More succinctly put, within the region $R$ we move $k$-tokens horizontally
from column $i$ to column $j$ and vertically from row $j$ to row $i$,
wherever possible.  See Figure~\ref{fig:flaggame}.

In the play of the game we may make any sequence
of moves in any order.  
The game is {\bf won} when there is exactly one token in each square.  

\begin{remark}
\label{rmk:maxsplit}
\rm
It turns out to be advantageous to split along an ideal subset $A$ 
if and only if the total number
of tokens in all squares in $A$ equals $|A|$.  We call 
the process
of finding all such $A$ and splitting along them {\bf splitting maximally}.
Although in this paper we won't take full advantage of this fact by 
splitting maximally before every move, we will never even
consider the possibility of splitting along $A$ if this condition is not 
met.
\end{remark}

The main result that we shall need is the following.

\begin{theorem}[\cite{P}]
\label{thm:win}
If the game for $\pi_1, \ldots, \pi_m$ can be won, then 
$$\int_{Fl(n)} [X_{\pi_1}] \cdots [X_{\pi_m}] \geq 1.$$
\end{theorem}

In general we do not know if the converse of Theorem~\ref{thm:win}
is true.  When $m=3$, the Schubert intersection numbers~\eqref{eqn:intnumber}
are structure constants for the ring $H^*(Fl(n))$, and in this case
the converse has been confirmed for $n \leq 7$.
It would certainly be remarkable if it were true in general.

We can also use the game to study the cohomology rings of partial
flag manifolds by pulling back cohomology classes to the full 
flag manifold.  In this paper, we investigate this in some detail in the
case of the Grassmannian.  Our main result is a version of the
converse of Theorem~\ref{thm:win} for Grassmannian Schubert calculus.

When $m=3$ and the classes come from a Grassmannian, the intersection
numbers~\eqref{eqn:intnumber} are Littlewood-Richardson numbers.
These numbers are also important in representation theory and
in the theory of symmetric functions---they 
are the structure constants for the representation ring of $GL(n)$, and 
for the ring of symmetric functions in the Schur basis.
As such, they are well studied,
and there are a number of combinatorial rules, and geometric rules
(e.g. \cite{C,V})
known for computing these numbers.

\begin{remark} \rm
An interesting and pleasant feature of root games is that the rules naturally
extend to any number Schubert classes in a way which is manifestly symmetrical
in these inputs: it is immediately clear from the definitions that reordering 
the input permutations $\pi_1, \ldots, \pi_m$ does not affect whether 
or not the game can be won.  This property is not generally
shared by other combinatorial
or geometric rules for Littlewood-Richardson numbers.  Most of the known 
rules have no
manifest symmetry.  Knutson-Tao puzzles \cite{KTW} are perhaps the
most manifestly symmetrical Littlewood-Richardson rule, having a $3$-cyclic 
symmetry when $m=3$, but
this is lost when one attempts to generalise beyond beyond triple 
intersections.   

The manifest symmetry of the root game on $Fl(n)$,
descends to the Grassmannian case.  In Section~\ref{sec:associating} 
we partially break this symmetry, but as we explain in 
Section~\ref{sec:remarks}, the full symmetry is easily
restored.  To our knowledge, the only other manifestly symmetrical rule 
in Grassmannian Schubert calculus is the Horn recursion (see the survey
article \cite{F2}), which, like the root game, does not explicitly compute
Littlewood-Richardson numbers and is only for determining which
Schubert intersection numbers are strictly positive.
\end{remark}

%
%%%%%%%%%%%%%%%%%%%%%%%%%%%%%%%%%%%%%%%%%%%%%%%%%%%%%%%%%%%%%%
%

\section{Associating a game to a Grassmannian Schubert calculus problem}
\label{sec:associating}

\begin{definition} \rm
A {\bf $01$-string} is a string 
$\sigma=\sigma^1 \ldots \sigma^n$  where each $\sigma^i \in \{0,1\}$.
A {\bf $0^m1^l$-string} is a $01$-string 
$\sigma=\sigma^1 \ldots \sigma^{m+l}$, where exactly $l$ of the $\sigma^i$
are equal to $1$.
\end{definition}

Schubert varieties in the Grassmannian 
$Gr_l(n)$ are indexed by $0^{n-l}1^l$-strings.  Fix a base flag
$$\{0\}=V_0 \subsetneq V_1 \subsetneq \cdots \subsetneq V_n =\CC^n$$ 
in $\CC^n$.
The Schubert variety $Y_\sigma \subset Gr_l(n)$
corresponding to $\sigma$ is 
$$Y_\sigma = 
\{y \subset\CC^n \ |\ \dim y \cap V_i \geq \sigma^1 + \cdots + \sigma^i\}.$$  
We denote its cohomology class in $H^*(Gr_l(n))$ by $[Y_{\sigma}]$.  
According to these
conventions $[Y_{0\ldots01\ldots1}]$ is the identity element in
$H^*(Gr_l(n))$ and $[Y_{1\ldots10\ldots0}] \in H^\text{top}(Gr_l(n))$
is the class of a point.

Given a list of $s+2$ $0^{n-l}1^l$-strings 
$\sigma_1, \ldots, \sigma_s, \mu, \nu$, we will wish to study the 
Grassmannian Schubert intersection numbers
$$\int_{Gr_l(n)} [Y_{\sigma_1}] \cdots [Y_{\sigma_s}] [Y_\mu] [Y_\nu].$$ 
We do so by investigating an equivalent problem on a full flag manifold.
The most obvious way to do this is to simply consider the product of
the classes $\alpha^*([Y_{\sigma_i}])$ etc., under the natural map 
$\alpha:Fl(n) \to Gr_l(n)$, as in Lemma~\ref{lem:pullbacka} below.  
However, our purposes require that we
do things in a somewhat less straightforward way.

Let $N \geq 0$ be an integer, and let $\sigma$ be a $0^{n-l}1^l$-string.  
Let
$i_1 < \cdots < i_{n-l}$ denote the positions of the zeroes in $\sigma$, and
$j_1 < \cdots < j_l$ denote the positions of the ones.
We define three ways to associate a permutation to the $01$-string $\sigma$:
\begin{align*}
\pi{(\sigma, N)}
&= i_1 \ldots i_{n-l} j_1 \ldots j_l (n{+}1) (n{+}2) \ldots (n{+}N) 
\\
\pi'{(\sigma, N)}
& = (i_1{+}N) \ldots (i_{n-l}{+}N) 1 2 \ldots N (j_1{+}N) \ldots (j_l{+}N)
\\
\pi''{(\sigma,N)} 
&= i_{n-l} \ldots i_1 (n{+}N) \ldots (n{+}1) j_l \ldots j_1.
\\
\end{align*}
From $\sigma_1, \ldots, \sigma_s, \mu, \nu$ we produce a 
list of permutations, $\pi_1, \ldots, \pi_{s+2} \in S_{n{+}N}$:
\begin{align*}
\pi_1 &= \pi{(\sigma_1,N)} \\
& \ \ \vdots \\
\pi_s & = \pi{(\sigma_s,N)} \\
\pi_{s+1} & = \pi'{(\mu,N)} \\
\pi_{s+2} &= \pi''{(\nu,N)}. \\
\end{align*}

\begin{proposition}
\label{prop:grtoflag}
$$\int_{Gr_l(n)} 
[Y_{\sigma_1}] \cdots [Y_{\sigma_s}] [Y_\mu] [Y_\nu] 
=
\int_{Fl(n{+}N)}
[X_{\pi_1}] \cdots [X_{\pi_{s+2}}].$$
\end{proposition}

The proof is based on the following standard pullback calculations,
whose proofs we omit.

\begin{lemma}
\label{lem:pullbacka}
Let $\alpha:Fl(n) \to Gr_l(n)$ be the map which forgets all but
the $l$-dimensional subspace of the flag.  Then
$\alpha^*([Y_\sigma]) = [X_{\pi(\sigma,0)}]$.
\end{lemma}

For a $01$-string $\sigma$, let $\sigma_+$ denote the string $\sigma$
followed by $N$ ones, and let $_+\sigma$ denote the string $\sigma$
preceded by $N$ ones.

\begin{lemma}
\label{lem:pullbackb}
Let $\beta:Gr_l(n) \to Gr_{l{+}N}(n{+}N)$ be the map 
$V \mapsto V \times \CC^N \subset \CC^n \times \CC^N$.
Then $\beta^*([Y_{\sigma_+}]) = [Y_\sigma]$.  If $\sigma'$ is not of 
the form $\sigma_+$ for some $0^{n-l}1^l$-string $\sigma$
then $\beta^*([Y_{\sigma'}]) = 0$.
\end{lemma}

If $\mathcal{F}$ is a partial flag variety and
$h \in H^*(\mathcal{F})$ is a Schubert class, let $h^\vee$ denote the 
opposite Schubert class, i.e. the unique Schubert class such that
$\int_\mathcal{F} h\cdot h^\vee = 1$.  For example 
$[X_\pi]^\vee = [X_{w_0\pi}]$, and 
$[Y_\sigma]^\vee = [Y_{\sigma^\text{rev}}]$ where 
$\sigma^\text{rev}=\sigma^n\ldots\sigma^1$ is $\sigma$ reversed.

If $h_1, \ldots, h_r \in H^*(\mathcal{F})$ are Schubert classes then the 
statement that
\begin{equation}
\label{eqn:unspecial}
\int_\mathcal{F} h_1 \cdots h_r = c 
\end{equation}
is equivalent to the statement that
\begin{equation}
\label{eqn:special}
h_1 \cdots \widehat{h_i} \cdots h_r = c\, h_i^\vee + \cdots
\end{equation}
in the Schubert basis.  We'll call Equation 
\eqref{eqn:special} the $h_i$-special
version of Equation \eqref{eqn:unspecial}.

\begin{proof}[Proof of Proposition~\ref{prop:grtoflag}]
Consider the equation in $H^*(Gr_{l{+}N}(n{+}N))$:
\begin{equation}
\label{eqn:middlestep}
\int_{Gr_{l{+}N}(n{+}N)}
[Y_{{\sigma_1}_+}] \cdots [Y_{{\sigma_s}_+}] [Y_{_+\mu}] [Y_{\nu_+}] = c.
\end{equation}

If we take the $[Y_{_+\mu}]$-special version of 
Equation \eqref{eqn:middlestep}
and pull it back to $Gr_l(n)$, we get (using Lemma~\ref{lem:pullbackb})
the $[Y_{\mu}]$-special version of 
$$\int_{Gr_l(n)} 
[Y_{\sigma_1}] \cdots [Y_{\sigma_s}] [Y_\mu] [Y_\nu] = c.$$

On the other hand, if we take the $[Y_{\nu+}]$-special version of 
Equation \eqref{eqn:middlestep} 
and pull it back to $Fl(n{+}N)$, we get (using Lemma~\ref{lem:pullbacka})
the $[X_{\pi_{s+2}}]$-special
version of 
$$\int_{Fl(n{+}N)} [X_{\pi_1}] \cdots [X_{\pi_{s+2}}] = c.$$
\end{proof}

%
%%%%%%%%%%%%%%%%%%%%%%%%%%%%%%%%%%%%%%%%%%%%%%%%%%%%%%%%%%%%%%
%

\section{Non-vanishing for Grassmannians}

Let $\sigma_1, \ldots, \sigma_s, \mu, \nu$ be $0^{n-l}1^l$-strings
(if $s=1$, we'll
write $\sigma$ instead of $\sigma_1$).  For any given $N \geq 0$, we
associate permutations $\pi_1, \ldots, \pi_{s+2}$ as before.  Our
goal in this section is to prove the following theorem.

\begin{theorem}
\label{thm:winconverse}
Take $N$ suitably large  ($N \geq l$ will always suffice).  The
root game corresponding to $\pi_1, \ldots, \pi_{s+2}$ can be won if and
only if 
$$\int_{Gr_l(n)} 
[Y_{\sigma_1}] \cdots [Y_{\sigma_s}] [Y_\mu] [Y_\nu]  \geq 1.$$
Moreover, only moves involving tokens labelled $1, \ldots, s$ are required.
\end{theorem}

\begin{remark}
\rm
The sufficient condition $N \geq l$ is not a sharp bound.  
This fact raises a number of interesting questions, which 
will be discussed in Sections~\ref{sec:remarks} and~\ref{sec:geometry}.
In the meantime the reader should not be alarmed by examples which
use smaller values of $N$.
\end{remark}

We shall first consider what happens in the case where $s=1$. 

Recall the correspondence between $01$-strings and Young diagrams.
Our Young diagrams will be in the French convention (the rows are left
justified and increase in length as we move down).  If $\sigma$ is a
$01$-string $\sigma$, let $r_i(\sigma)$ denote the number of ones
before the $i^{\rm th}$ zero.  We associate to $\sigma$ the Young
diagram $\lambda(\sigma)$ whose $i^{\rm th}$ row is $r_i(\sigma)$,
where we are allowing the possibility that some rows may have length 
$0$.

If $\lambda$ is a Young diagram, let $N+ \lambda$ denote the
Young diagram obtained by adding $N$ squares to each row of $\lambda$
{\em including those rows which contain $0$ squares}.

\begin{figure}[htb]
  \begin{center}
    \epsfig{file=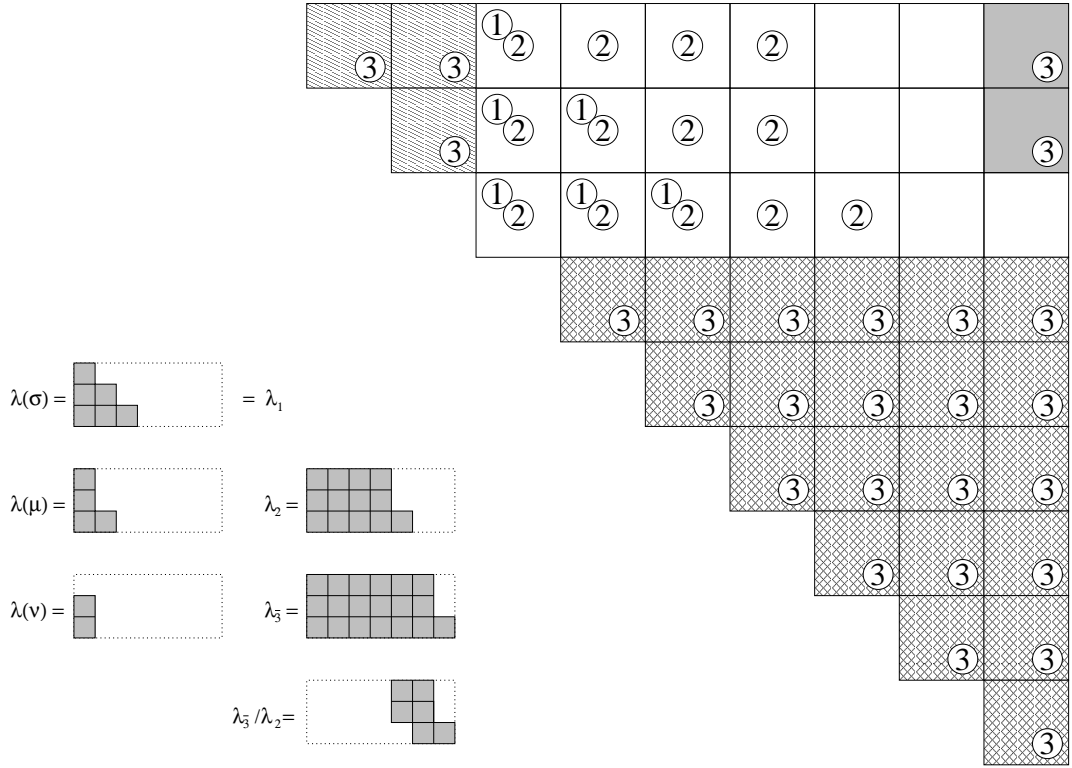,height=4in}
    \caption{Initial position of the game for $\sigma=1010101$, 
        $\mu=1001011$, $\nu=0100111$, with $N=3$.  Here
        $\pi_1 = 246135789\ten$, $\pi_2 = 568123479\ten$, 
        $\pi_3 = 431\ten987652$.}
    \label{fig:grinitial}
  \end{center}
\end{figure}

The initial positions of the $1$-tokens are in the shape
of the Young diagram $\lambda_1 = \lambda(\sigma)$.  The initial positions
of the $2$-tokens are in the shape of a Young diagram 
$\lambda_2 = N+\lambda(\mu)$.  The squares that do {\em not} contain
a $3$-token are also in the shape of a Young diagram, which we'll denote
$\lambda_{\bar 3}$; viewed upside down, $\lambda_{\bar 3}$ is the 
complement to $\lambda(\nu)$ inside
an $(n-l) \times (l+N)$ rectangle.  The lower left corner of
each of the Young diagrams $\lambda_1$, $\lambda_2$, $\lambda_{\bar 3}$
is in the square $S_{n-l, n-l+1}$.
See Figure~\ref{fig:grinitial} for an illustration
of how these shapes are generated.

If $\lambda_2 \nsubseteq \lambda_{\bar 3}$, then it is a basic fact that
$[Y_{\mu}] [Y_{\nu}] = 0 \in H^*(Gr_l(n))$ (in fact this is a necessary
and sufficient condition).  Therefore, we may assume that
no square contains both a $2$-token and a $3$-token.
The squares which contain neither a $2$-token
nor a $3$-token are empty squares---since $N$ is suitably 
large the $1$-tokens are all to the left of these squares---
and are in the shape of a skew-diagram $\skewtabl$.

At the outset of the game, some immediate 
splitting can occur.  We split in such a way that each square containing 
a $3$-token becomes a $1$-square
region of its own.  (For some choices of $(\sigma, \mu, \nu)$ it may
be possible to split beyond this, but our argument is slightly 
simplified if we elect not to.)
The remaining squares are those of $\lambda_{\bar 3}$,
which form what we call the {\bf big region}.  
The big region is the only region which is unsolved; naturally, therefore,
this will be the region in which all moves take place.

It is worth taking a moment to note how tokens move within the
big region.  A priori, a move $(i,j)$ will cause some
$k$-tokens to move horizontally and others to move vertically.
However, since the rows of the big region
are indexed by $\{1, \ldots, n-l\}$, and the columns are 
indexed by the disjoint set $\{n-l+1, \ldots, n+N\}$,
these cannot both happen.
If $i<j\leq n-l$ then $k$-tokens
will move vertically from row $j$ to row $i$.  
If $ n-l+1 \leq i < j$ then 
$k$-tokens move horizontally from column $i$ to column $j$.
No tokens move if $i \leq n-l < j$.

\begin{definition}[Zelevinsky \cite{Z}]  
\label{def:picture} \rm
 A {\bf picture} between two (French) skew diagrams is a bijection
between their boxes with the property that if box $A$ is weakly
above and weakly right of box $B$ in one diagram, then the corresponding
boxes $A'$ and $B'$ are in lexicographic order (i.e. $S_{ij}$ precedes  
$S_{i'j'}$ if $i<i'$ or $i=i'$ and $j<j'$) in the other diagram.
\end{definition}

Now
$\int_{Gr_l(n)} [Y_\sigma] [Y_\mu] [Y_\nu]$
is given by the Littlewood-Richardson coefficient 
$c_{\lambda_1 \lambda_2}^{\lambda_{\bar 3}}$, which
can be described as the number of pictures between 
$\lambda_1$ and $\skewtabl$ \cite{Z}.
This is a reformulation of the Littlewood-Richardson rule \cite{LR},
closely related to the Remmel-Whitney formula \cite{RW}.
In particular, if this number is non-zero, there exists such a picture. We
pick one, and denote by $f$ the map it defines from the squares of 
$\lambda_1$ to the squares of $\skewtabl$.
Note that all of these squares are in the big region of the game.  

To show that the game can be won, we will give an algorithm---
the {\em Grassmannian root game algorithm} (GRGA)---which uses 
$f$ to construct a sequence of moves that transports each $1$-token 
to a square of $\skewtabl$.

Essential to the GRGA is the following numbering scheme.
At each point in the game we associate a number---called the
{\em readiness number}---to each unplaced $1$-token 
(i.e. one which has not already 
reached its final destination) and each empty square of $\skewtabl$.

\begin{definition}
\rm
At any point in the game, let $t$ be a $1$-token whose initial square 
was $S \in \lambda_1$
and whose current square is $S_{ij}$.  
Let $S_{i'j'} = f(S)$.
Define the {\bf readiness number} of 
both the token $t$ and the square $S_{i'j'}$ to be the number $i-i'$.
We say an empty square of $\skewtabl$ or an unplaced $1$-token 
is {\bf ready} if its readiness number is $0$.  Tokens which
have reached their final destination and non-empty squares are
{\em not} considered ready.
\end{definition}

The key properties of this numbering scheme are the following:

\begin{lemma}
\label{lem:numbering}
Initially, the readiness numbers of the unplaced $1$-tokens are 
\begin{enumerate}
\item[(a)] weakly increasing along each row, and  
\item[(b)] weakly decreasing down each column;
\end{enumerate}
in $\skewtabl$ the readiness numbers of empty squares are 
\begin{enumerate}
\item[(c)] weakly decreasing along each row, and 
\item[(d)] weakly increasing down each column. 
\end{enumerate}
Moreover, the GRGA, described below, preserves all of these properties.
\end{lemma}

\begin{proof}
(a) If $A$ and $B$ are squares of $\lambda_1$ in the same row, and
$A$ is right of $B$, then by Definition~\ref{def:picture} $f(A)$ is
lexicographically before $f(B)$.  In particular, $f(A)$ is weakly
above $f(B)$, i.e. in the same row or a row above.  Thus the readiness
number of the token in $B$ $\leq$ the readiness number of the token
in $A$.

(b) If $A$ and $B$ are squares of $\lambda_1$, and $A$ is one square
above $B$, then again $f(A)$ is lexicographically before $f(B)$.
There are two cases. If $f(A)$ is strictly above $f(B)$, then the
readiness number of the token in $B$ $\leq$ the readiness number of 
the token in $A$.  Otherwise,
$f(A)$ and $f(B)$ are in the same row in $\skewtabl$, with $f(B)$
right of $f(A)$.  But then by Definition~\ref{def:picture}, $B$ 
must be lexicographically before $A$, which is a contradiction.

Statements (c) and (d) are proved similarly.  
That the GRGA preserves all these properties will be quite evident.
\end{proof}

Note that since the readiness number of the lower-leftmost token is at
least $0$, by Lemma~\ref{lem:numbering} parts~(a) and~(b) 
the readiness numbers are initially all non-negative. 

\paragraph{The Grassmannian root game algorithm (GRGA).}
{\em
The algorithm assumes that the tokens are in the initial positions of
the root game for $\pi_1, \pi_2, \pi_3$ (corresponding to 
$\sigma, \mu, \nu$ with $N$ suitably large), that
all $3$-tokens have been split into their own one-square region,
and that we have a picture $f$ between $\lambda_1$ and $\skewtabl$.
All moves take place in the big region.}

\begin{enumerate}
\item If any of the $1$-tokens are ready, go to Step 2.
Otherwise, perform a sequence of moves to 
shift all unplaced $1$-tokens up one square.  
The reader can easily check that the sequence of moves 
$(1,2),\ (2,3),\ \ldots,\ (n-l-1,n-l)$ accomplishes this.
The assumption that $N$ is sufficiently
large ensures that the upward movement of the $1$-tokens
is unobstructed.
This step will decrease the readiness number of each $1$-token by $1$.  
Repeat this step until some $1$-token is ready.

\item Scan through the columns of $\skewtabl$, beginning with the
rightmost column and proceeding to the left.  Within each column locate
the topmost square that does not already contain a $1$-token.  Let $S$ be
the {\em first} ready square which we encounter in this way.

\item Find a ready token $t$ in the same row as the square $S$.  Make
the unique move which causes $t$ to move into $S$.  This may cause other
tokens to move as well.  All tokens which move here move to 
their final destination, so after this move, they and the squares
they occupy are no longer considered ready.

\item Repeat Steps 1 through 3 until every square of 
$\skewtabl$ contains a $1$-token.
\end{enumerate}

%For $s=1$ the game will be won at this point.  However in order to make
%use of this algorithm effectively in the case where $s>1$, we will
%require the following additional step.
%\begin{enumerate}
%\item[5.] Split in a minimal way so that every $1$-token is in a 
%one-square region of its own.
%\end{enumerate}
%
%A small example of the GRGA and all of its the moves from beginning
%to end is illustrated in Example \ref{ex:algorithm}, at the end of this 
%section.

\begin{example}
\label{ex:algorithm}
\rm
Figures~\ref{fig:algorithm} and~\ref{fig:algorithmcont} illustrate
the GRGA, applied to the example from Figure~\ref{fig:grinitial}. 
We now draw only the squares in the upper right $3 \times 7$
rectangle as these are the only ones relevant to the movement of the 
$1$-tokens.
Moreover, only the $1$-tokens are shown in these diagrams, and the number
on the token is the readiness number, not the token label.  To specify
the picture $f$, each 
$1$-token is given a shading and the corresponding square under $f$ in 
$\skewtabl$ is shaded similarly.  Each unshaded square actually contains a 
$2$-token.  The two darkly shaded squares in the upper right corner contain 
$3$-tokens, as do each of the squares not shown in this diagram, but
these squares are not part of the big region.
\end{example}

\begin{figure}[htbp]
  \begin{center}
    \epsfig{file=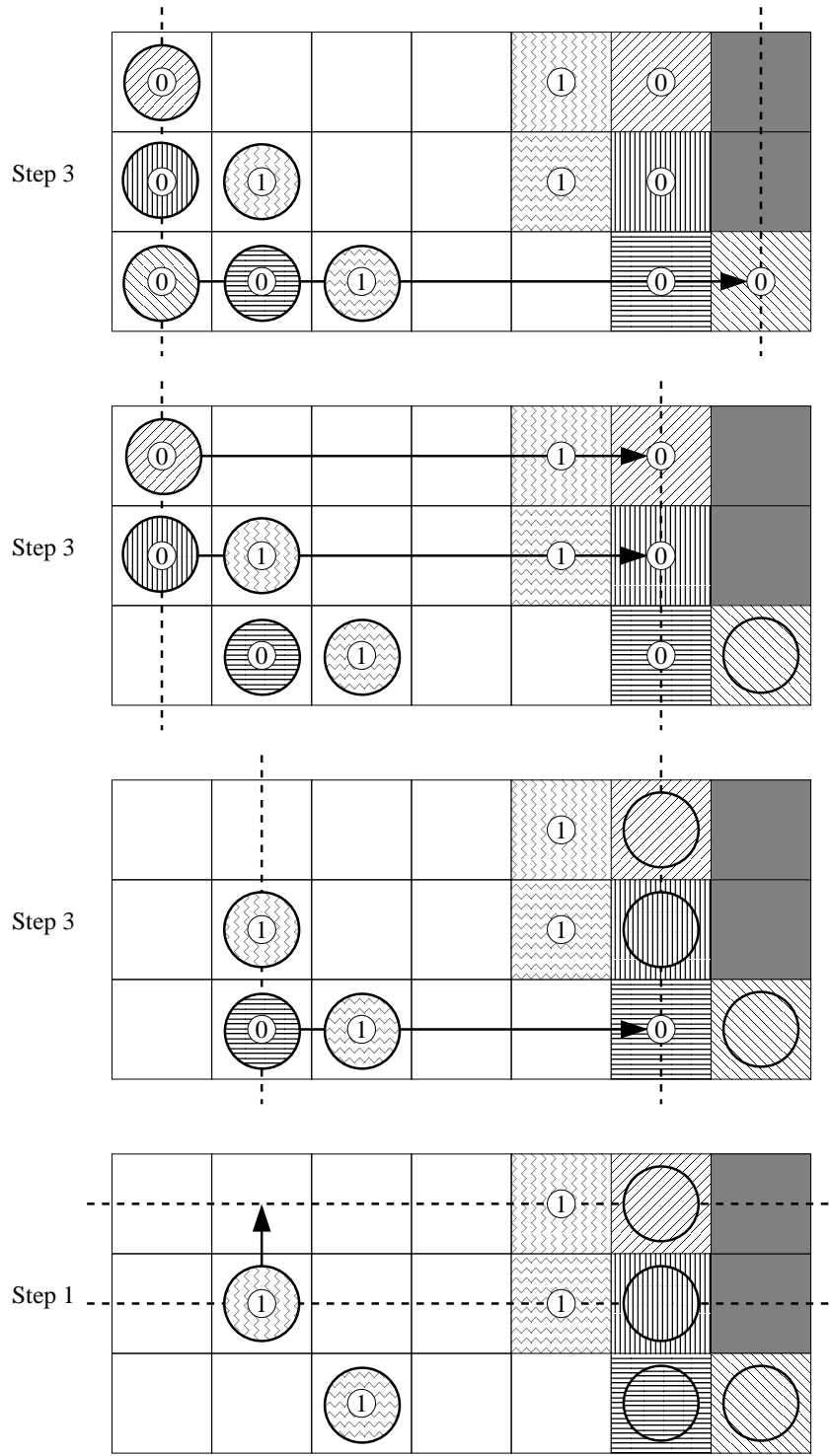,height=7.7in}
    \caption{The Grassmannian root game algorithm.  Here
	$\sigma=1010101$, $\mu=1001011$, $\nu=0100111$, and $N=3$.}
    \label{fig:algorithm}
  \end{center}
\end{figure}
\begin{figure}[htbp]
  \begin{center}
    \epsfig{file=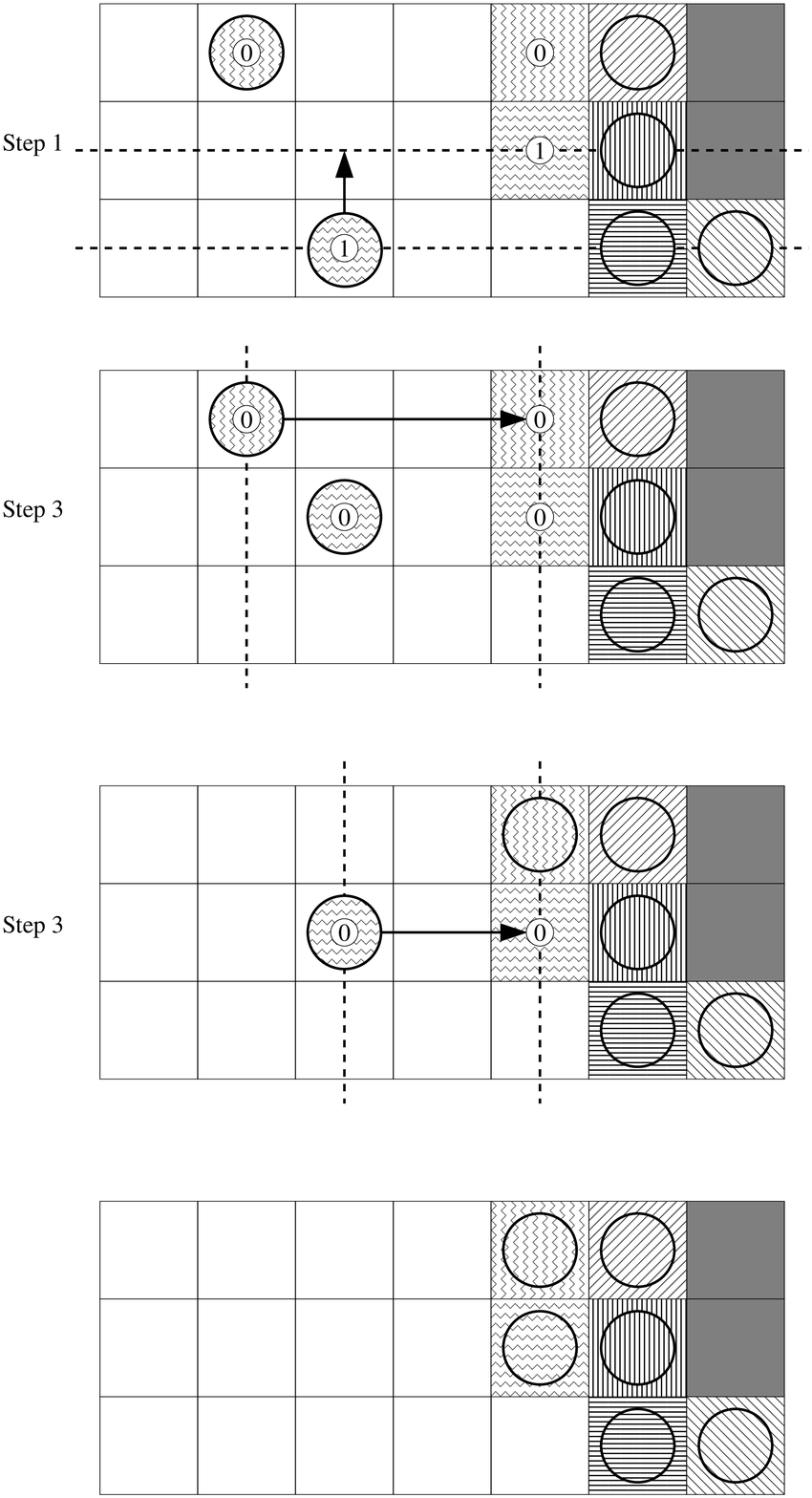,height=7.7in}
    \caption{Continuation of Figure~\ref{fig:algorithm}}
    \label{fig:algorithmcont}
  \end{center}
\end{figure}

We now show that the GRGA accomplishes what it claims to accomplish.

\begin{lemma}
\label{lem:3converse}
Given a picture $f$ between $\lambda_1$ and $\skewtabl$, the GRGA
will win the game for $\pi_1, \pi_2, \pi_3$.
\end{lemma}

\begin{proof}
First note that after vertical movement from Step 1 is finished, the
number of ready tokens in any row equals the number of ready squares in
that row: an empty square $S$ is ready if and only if the $1$-token which 
began in the square $f^{-1}(S)$ is in the same row as $S$.  
We show that this equality is preserved, by showing that
the move in Step 3 only ever causes ready tokens to move into ready
squares.

\paragraph{Claim (i)} 
{\em A move from Step 3 causes only ready tokens to move.
In particular the number of ready tokens in a row always remains less 
than or equal to the number of ready squares in a row.  }
The only tokens that can conceivably move are those in the same column 
as $t$.
Because $S$ is the top unfilled square in its column, no tokens above 
$t$ move.  
Because of Lemma~\ref{lem:numbering} part~(b), 
all tokens below $t$ are ready.

\paragraph{Claim (ii)} 
{\em Only ready squares are filled.} 
The algorithm attempts to fill the rightmost squares first.  If there is
a ready square $S'$ in some column, the topmost empty square 
in that column will also be ready, by  Lemma~\ref{lem:numbering} part~(d);
thus the algorithm will never fill any square left of $S'$ before it
fills $S'$. 
However, by Lemma~\ref{lem:numbering} part~(c) the ready squares 
are rightmost in
their row.  Thus if a token moves into a non-ready square, it means
that there are no ready squares in its row.  But since only ready
tokens move, we would have a row with at least one ready token
and no ready squares.  This, as noted in Claim~(i), is impossible.

\paragraph{Claim (iii)} {\em The move from Step 3 is always possible.}
Since only ready squares are filled by ready tokens, the number of 
ready squares and ready tokens in any given row is always equal.  Thus
there is a ready token $t$ in the same row as the ready square $S$.
Because we assume $N$ to be sufficiently
large, $t$ is to the left of $S$.  If $t$ is in column $i$ and
$S$ is in column $j$, the move $(i,j)$ will take the token $t$ 
into square $S$.

Thus in Steps~2 and~3, every ready square eventually gets filled by a 
ready token:
by the argument in Claim~(ii) no square is skipped.
However, because the readiness numbers are initially non-negative and
Step~1 decreases the readiness number of each square 
by $1$, every square of $\skewtabl$ is ready at some point; 
thus the algorithm puts a $1$-token in each square of $\skewtabl$, at
which point the game is won.
\end{proof}

\begin{proof}[Proof of Theorem~\ref{thm:winconverse}]
($\,\Longrightarrow\,$)
This follows from 
Proposition~\ref{prop:grtoflag} and Theorem~\ref{thm:win}.

($\,\Longleftarrow\,$)
For $s=1$ we use the GRGA, which wins the game by Lemma~\ref{lem:3converse}.
For $s>1$, we proceed by induction.
Suppose $\int_{Fl(n{+}N)} [X_{\pi_1}] \cdots [X_{\pi_{s+2}}] \neq 0$.
Then we can write
\begin{equation}
\label{eqn:induction}
[X_{\pi_2}] \cdots [X_{\pi_{s+1}}] = c [X_{\rho}] + \cdots 
\end{equation}
in the Schubert basis, where $c>0$, and
$$\int_{Fl(n{+}N)} [X_{\pi_1}] [X_{\rho}] [X_{\pi_{s+2}}] \neq 0.$$
Since this is really a Grassmannian calculation, $\rho$ will be
necessarily be of the form 
$\pi'(\sigma',N)$ for some $0^{n-l}1^l$-string $\sigma'$.
By Lemma~\ref{lem:3converse} we can win the game corresponding to 
$\pi_1, \rho, \pi_{s+2}$, only moving $1$-tokens.
It is easy to see that exactly the 
same sequence of splittings and moves can be made in the 
game for $\pi_1, \pi_2, \ldots, \pi_{s+1}, \pi_{s+2}$,
and that it causes the $1$-tokens to end up in exactly the same final 
positions.  Note that we end up with either a $1$-token or an 
$(s{+}2)$-token in every
square which does not correspond to an inversion of $\rho$, i.e. every
square which does correspond to an inversion of $w_0 \rho$.

This sequence of moves no longer wins the game; however, we can proceed
inductively, after two further small steps.  First, we perform a sequence 
of splittings so that every $1$-token is in a one-square region of its
own.  Next, after splitting in this way, we replace each
$1$-token by an $(s{+}2)$-token.  This second step is not a legitimate play
in the game, but it is completely harmless: 
because every $1$-token is sequestered in its own 
one-square region, 
it can have no effect whatsoever 
on any possible subsequent moves of the game. 
But now we have precisely reached the initial 
position of the game 
corresponding to $\pi_2, \ldots, \pi_{s+1}, w_0 \rho$.
This is again a game associated to a Grassmannian problem, and by Equation
\eqref{eqn:induction} the Schubert intersection number is non-zero.  
By induction, there is a sequence of moves to win this new game.  Thus
by concatenating the two sequences of moves, we can win the original game.
\end{proof}

%
%%%%%%%%%%%%%%%%%%%%%%%%%%%%%%%%%%%%%%%%%%%%%%%%%%%%%%%%%%%%%%
%

\section{Remarks}
\label{sec:remarks}
In Step 3 of the GRGA, there is a somewhat canonical choice for the
token $t$, namely the leftmost ready token in its row.  If we use this 
choice of $t$, one can verify that the algorithm actually transports
the $1$-token which is initially in square $S$ to the square $f(S)$.
On the other hand, if there is a way of winning the root game, there
is generally a plethora of ways, most of which do not arise by 
following the GRGA for any picture.  
Theorem~\ref{thm:winconverse} tells us that the existence
of any one way to win implies the existence of a picture between
$\lambda_1$, and $\skewtabl$.  However, given a sequence of moves
which wins the game, it is not at all obvious how to construct such
a picture.  This is even unclear in the the case where $s=1$, and
only $1$-tokens are moved.

It is worth noting that the root game can be used to determine whether
$$[Y_{\sigma_1}] \cdots [Y_{\sigma_s}] [Y_\mu] [Y_\nu] \neq 0$$
even if the cohomological degree of the product is not $\dim_\RR Gr_l(n)$.
To do this, we modify the game by changing the winning condition to
read ``the game is won if there is {\em at most} one token in each square'', 
rather than ``exactly one token in each square''.  Once we do this, we 
have the following corollary of Theorem~\ref{thm:winconverse}.
\begin{corollary}
Take $N$ suitably large, and let 
$\pi_1, \ldots, \pi_{s+2}$ be obtained from
$\sigma_1, \ldots, \sigma_s,\linebreak[0]\mu, \nu$ as before.  Then
$[Y_{\sigma_1}] \cdots [Y_{\sigma_s}] [Y_\mu] [Y_\nu] \neq 0$ if
and only if the root game for $\pi_1, \ldots, \pi_{s+2}$ can
be won with the modified winning condition.
\end{corollary}
\begin{proof}
Assume
$[Y_{\sigma_1}] \cdots [Y_{\sigma_s}] [Y_\mu] [Y_\nu] \neq 0$. Then 
there exists $\sigma'$ such that $$\int_{Gr_l(n)} 
[Y_{\sigma'}] [Y_{\sigma_1}] \cdots [Y_{\sigma_s}] [Y_\mu] [Y_\nu] \neq 0.$$
Let $\pi' = \pi(\sigma',N)$.  Since we can win the unmodified game for
$\pi', \pi_1, \ldots, \pi_{s+2}$, we can win the modified game
for $\pi_1, \ldots, \pi_{s+2}$ simply by omitting 
moves where the token corresponds to $\pi'$.  
The reverse direction follows from Proposition~\ref{prop:grtoflag} 
and \cite[Theorem 5]{P} (which generalises Theorem~\ref{thm:win}).
\end{proof}
\noindent
There is a small catch: with this more general winning condition, our 
observation in Remark~\ref{rmk:maxsplit} becomes invalid.  There is no 
longer an easy necessary and sufficient condition indicating when splitting 
is advantageous.

\begin{figure}[htb]
  \begin{center}
    \epsfig{file=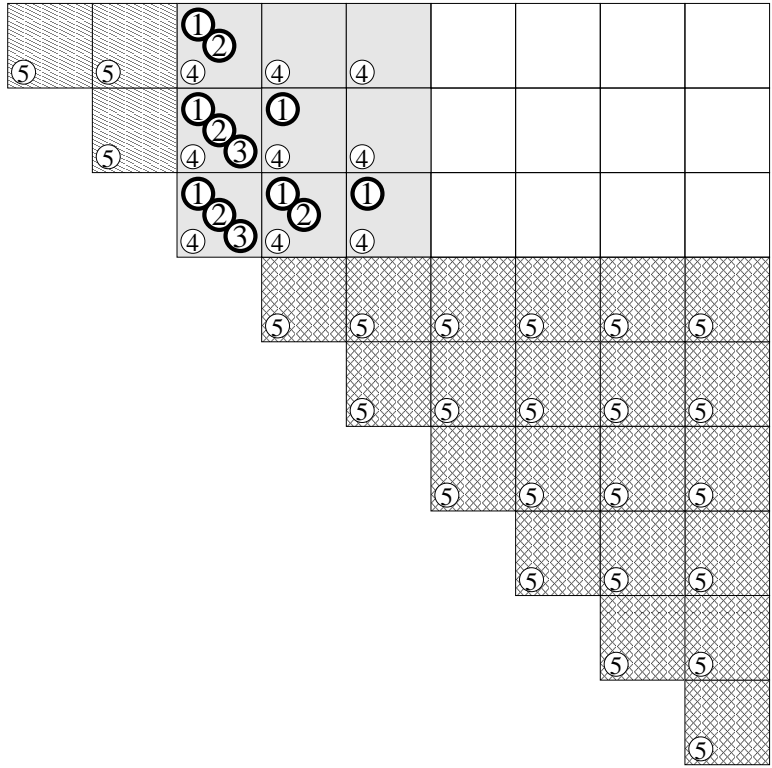,height=4in}
    \caption{Initial position of the game for $\sigma_1=1010101$, 
	$\sigma_2=1001011$, $\sigma_3=0100111$, $\mu=\nu=0001111$, with $N=3$.
	Squares are shaded if they contain a $4$-token or a $5$-token. 
	This is the symmetrical version of the example in
	Figure~\ref{fig:grinitial}.}
    \label{fig:grinitialsymm}
  \end{center}
\end{figure}

One of the unfortunate features of this presentation is the asymmetry
in the way the permutations $\pi_1, \ldots, \pi_{s+2}$ are defined.
The root game itself is manifestly symmetrical in the permutations
given.  However, because $\pi_{s+1}$ and $\pi_{s+2}$ are produced
in a different way from $\pi_1, \ldots, \pi_s$, the symmetry
is broken for Grassmannians.
Nevertheless, as Theorem~\ref{thm:winconverse} is valid for any
$s$, we can formulate a 
symmetrical game by taking $\sigma_1,\ldots,\sigma_s$ to be arbitrary, and
$\nu = \mu = 0 \ldots 01 \ldots 1$, so that 
$[Y_\nu] = [Y_\mu] = 1 \in H^*(Gr_l(n))$.
To see how this new formulation changes the initial position,
contrast Figure~\ref{fig:grinitialsymm} 
with Figure~\ref{fig:grinitial}.

The only splitting which occurs in the proof of Theorem~\ref{thm:winconverse} 
is before the first move, and in the inductive
step.  The GRGA itself does not split between moves.
As noted in Remark~\ref{rmk:maxsplit}, it can never
be harmful to split maximally between moves, and it turns out that
if one modifies the GRGA to split maximally between moves, things
proceed very much as before.  However, in the next section, our proof of
Theorem~\ref{thm:geomconverse} will rely on the fact that the
GRGA involves no splitting.

It would be nice if we could take $N=0$ in Theorem~\ref{thm:winconverse}.
Although we are not aware of any example which proves that this cannot
be done, the algorithm simply falls apart if $N$ is too small.  There
are several problems which occur with trying to follow a similar
approach. The most serious of these is that a token may be to the
right of the square for which it is supposedly destined
according to the chosen picture.  Again, this highlights the fact
that we do not know a straightforward two-way correspondence between
pictures and ways of winning the root-game.
Instead, in the next section we prove 
Theorem~\ref{thm:geomconverse},
a geometrical analogue of Theorem~\ref{thm:winconverse}
which is valid for all $N$.  Theorem~\ref{thm:geomconverse} suggests
that it is not unreasonable to conjecture that 
Theorem~\ref{thm:winconverse} is true for all $N$.  We leave it
as an open problem to determine whether or not this is in fact the
case.

%
%%%%%%%%%%%%%%%%%%%%%%%%%%%%%%%%%%%%%%%%%%%%%%%%%%%%%%%%%%%%%%
%

\section{Geometric interpretation}
\label{sec:geometry}
In \cite{P} we give a complete description of the geometry underlying the 
root game.  
The picture is quite different from those found in the geometric 
Littlewood-Richardson rules of 
Vakil~\cite{V} and more recently Coskun~\cite{C}, both of which study 
degenerations of intersections of Schubert varieties
inside a Grassmannian---we would be surprised if there were any
straightforward relationship.  Our methods are based on studying
tangent spaces to Schubert varieties, and are more closely 
related to the
approach used by Belkale in his geometric proof of the Horn conjecture
\cite{B}.
Here, we will recall only the parts of the picture which are most 
relevant to our situation.

In this section we shall once again assume $s=1$.
Our notation changes slightly from Section~\ref{sec:prior} in that
we will be working with $GL(n{+}N)$ instead of $GL(n)$.  

Let $x_1, \ldots, x_{n{+}N}$ denote
the standard basis of $\CC^{n{+}N}$.
Let $B$ denote the standard Borel subgroup of $GL(n{+}N)$ (upper triangular
matrices), and let $B_-$ denote its opposite (lower triangular matrices).
As before $T = B \cap B_-$ will be the standard maximal torus.
For any complex vector space 
$V$ let $Gr(V)$ be the disjoint union of all Grassmannians 
$Gr_d(V)$, $0 \leq d \leq \dim V$.

Let $R$ be a region in the game, and let 
$$\tau:R \to \text{subsets of $\{1, 2, 3\}$}$$
describe the configuration of the tokens within this region.  In
the underlying geometry there is, assigned 
to the combinatorial pair $(R,\tau)$, a corresponding geometric 
pair $(V,U)$, where $V$ is a $B$-module, and $U=(U_1,U_2,U_3)$ is a 
$T^3$-fixed point on $Gr(V)^3$ (or equivalently the $U_k$ are 
$T$-invariant subspaces of $V$).  As a $T$-representation, $V$ is 
multiplicity-free, and the distinct $T$-weights 
are $$\text{weights}(V) = \{x_i - x_j\ |\ S_{ij} \in R\}.$$
The (distinct) $T$-weights
of $U_k$ correspond to the positions of the $k$-tokens inside $R$:
$$\text{weights}(U_k) = \{x_i - x_j\ |\ k \in S_{ij}\}.$$
Thus the pair $(V,U)$ carries all relevant information about the
region $R$ and the arrangement of the tokens with $R$.

The region $R$ is solved when there is exactly one token in
each square.  In terms of the pair $(V,U)$ this is occurs when
\begin{equation}
\label{eqn:transversedef}
V = U_1 \oplus U_2 \oplus U_3 .
\end{equation}
We'll call any $U = (U_1, U_2, U_3)$ which satisfies 
condition \eqref{eqn:transversedef} {\bf transverse}.

Assuming we do not split the region $R$, a move or a sequence of
moves in the game takes the pair
$(V,U)$ to a new pair $(V,U')$,
{\em where $U'$ is in the $B^3$-orbit closure through $U \in Gr(V)^3$}.
%If only $k$-tokens are moved then only $U'_k$ differs from $U_k$.
Thus if we solve a region starting from position $(V,U)$, we have
located a transverse $T^3$-fixed point $U' \in \overline{B^3 \cdot U}$.  

The importance of transverse points in $\overline{B^3 \cdot U}$
is seen in the following proposition.

%\begin{fact}
%\label{fact:fixedpt}
%For each region of the game consider the associated pair $(U,V)$.
%If for each region in the game there exists $T$-fixed point 
%$U' \in \overline{B^3 \cdot U}$ which is transverse, then
%$$\int_{Fl(n)} [X_{\pi_1}] [X_{\pi_2}] [X_{\pi_3}] \geq 1$$
%\end{fact}
%
%This is true at any stage in the game.  However from here on in we shall
%restrict our attention to the state of the game before the first
%move is made, but after splitting.  In which case, Fact 
%\ref{fact:fixedpt} is just one direction of the following result.

\begin{proposition}[\cite{P}]
\label{prop:transpt}
Consider the position of the root game 
game for $\pi_1, \pi_2, \pi_3$ which arises after
splitting but before the first move is made.  The tokens are in
their initial position, but there may be more than one region.
To each region $R$ there is an associated pair $(V_R,U_R)$.  Then
$$\int_{Fl(n)} [X_{\pi_1}] [X_{\pi_2}] [X_{\pi_3}] \geq 1$$
if and only if for every region $R$ there exists a transverse
point $U' \in \overline{B^3 \cdot U_R}$.
\end{proposition}

%The details of both of these facts are spelled out in the proof of 
%\cite[Theorem 3]{P} including how the pair $(U,V)$ is associated to an
%arrangement of tokens in a region.

Note that the point $U'$ in Proposition~\ref{prop:transpt} is not 
necessarily $T^3$-fixed.  The big question, therefore, is how specialised
can we make the point $U'$ and still have Proposition~\ref{prop:transpt} 
be true.  
There are three levels of specialisation that we could request
of this transverse point $U'$.

\begin{enumerate}
\item $U'$ is any transverse point in $\overline{B^3 \cdot U_R}$.
\item $U'$ is a $T^3$-fixed transverse point in $\overline{B^3 \cdot U_R}$.
\item $U'$ is a ($T^3$-fixed) transverse point in $\overline{B^3 \cdot U_R}$, 
where $(V_R,U')$ comes from applying sequence of root game moves (but no
splitting) to the position $(V_R,U_R)$.
\end{enumerate}

A priori, it is not clear that these three levels of specialisation
are equivalent.  However, for Grassmannian Schubert calculus with $N$ 
sufficiently large, Theorem~\ref{thm:winconverse} shows that they
are all equivalent.  The GRGA tells us exactly how to
produce a sequence of moves which gives the point $U'$ at
Level~3, which is the most specialised.  Note it is important here that
splitting is never used in the GRGA---when a region is split, the $U'$ 
one is tempted to define need not be in $\overline{B^3 \cdot U_R}$.

Unfortunately, when $N$ is small the GRGA can fail, and so we 
cannot claim that all three levels of specialisation are equivalent
for all $N$.  Our goal in this section is to show that even if $N$ is
too small for the GRGA to work, we can still get $U'$ at Level~2;
i.e. Proposition~\ref{prop:transpt} is still true for 
Grassmannian Schubert calculus if we demand that $U'$ be a $T^3$-fixed point.

%By Theorem \ref{thm:winconverse}, we know this is possible if
%$\int_{Gr_l(n)} [Y_{\sigma}] [Y_{\mu}] [Y_{\nu}] \geq 1.$
%
%Level 2 is equivalent to asking for a converse to Fact \ref{fact:fixedpt}.
%Since level 2 is less specialised than level 3, this converse is also true
%for Grassmannian Schubert problems if $N$ is sufficiently large.
%Our goal in this section is to show the that the converse of 
%Fact \ref{fact:fixedpt} is true in the case of Grassmannian Schubert 
%problems even for $N=0$. 

To make matters more concrete, we now explicitly describe the 
initial pair $(V,U)$ for the big region $R$ in the root game associated 
to $\sigma, \mu, \nu$.  This is the only region that we need to
concern ourselves with, since all other regions are already solved.

Let $M(n{+}N)$ be the space of $(n{+}N) \times (n{+}N)$ matrices, having 
standard basis $\{\tilde e_{ij}\}$, 
and let $B$ act on $M(n{+}N)$ by conjugation.
Let $W$ denote the $B$-submodule of $M(n{+}N)$ 
generated by the entries in the upper right $(n-l) \times (l+N)$ 
rectangle.  Let $W'$ be the $B$-submodule of $W$ generated
by $\tilde e_{ij}$ such that $S_{ij}$ contains a $3$-token.  Then
$V$ is the quotient $B$-module
$$V = W/W'.$$
Note $V$ has a basis 
$\{e_{ij}:= \tilde e_{ij}+W'\ |\ S_{ij} \in R\}$.  The point 
$U \in Gr(V)^3$ is described as follows:  
$$U_k = \text{span}\{e_{ij}\ |\ k \in S_{ij}\}.$$
Note that $U_3 = \{0\}$, so we need not give it much further consideration.

Let $\Udest$
be the subspace of $V$ whose $T$-weights correspond to $\skewtabl$:
$$\Udest= \text{span}\{e_{ij}\ |\ 2 \notin S_{ij}\}$$
and let $U' = (\Udest, U_2, U_3)$.
Note that $V = \Udest\oplus U_2 \oplus U_3$, so the point
$U' \in Gr(V)^3$ is transverse.

\begin{example} \rm
For the initial position shown in Figure~\ref{fig:grinitial},
\begin{align*}
V 
&=
\text{span} \{
\tilde e_{ij}\ |\ 1 \leq i \leq 3,\ 4 \leq j \leq 10\}
\,\big/\,\text{span} \{ 
\tilde e_{1\,\scriptten}, \tilde e_{2\,\scriptten}\} \\
&= 
\text{span} \{
e_{14}, e_{15}, \ldots, e_{19}, \:
e_{24}, e_{25}, \ldots, e_{29}, \:
e_{34}, e_{35}, \ldots, e_{39}, e_{3\,\scriptten}
\}, \\
U_1 &= \text{span} \{
e_{14}, \:
e_{24}, e_{25}, \:
e_{34}, e_{35}, e_{36} \}, \\
U_2 &= \text{span} \{
e_{14}, e_{15}, e_{16}, e_{17},\:
e_{24}, e_{25}, e_{26}, e_{27}, \:
e_{34}, e_{35}, e_{36}, e_{37}, e_{38} \}, \\
\Udest &= \text{span} \{
e_{18}, e_{19}, \:
e_{28}, e_{29}, \:
e_{39}, e_{3\,\scriptten}
\}.
\end{align*}
\end{example}

\begin{theorem}
\label{thm:geomconverse}
For every $N \geq 0$,
$$\int_{Gr_l(n)} [Y_{\sigma}] [Y_{\mu}] [Y_{\nu}] \geq 1$$
if and only if with $(V,U)$ and $\Udest$ as above, 
$\Udest \in \overline{B \cdot U_1}$.
\end{theorem}

\begin{proof}
$(\,\Longleftarrow\,)$ 
If $\Udest \in \overline{B \cdot U_1}$ then $U'$ is transverse point 
in $\overline{B^3 \cdot U}$.  Hence this follows from
Propositions~\ref{prop:grtoflag} and~\ref{prop:transpt}.

$(\,\Longrightarrow\,)$ 
Assume $\int_{Gr_l(n)} [Y_{\sigma}] [Y_{\mu}] [Y_{\nu}] \geq 1.$

We know the result is true for $N$ sufficiently large, since
the GRGA tells us how to get from the position $(V,U)$ to the
position $(V,U')$.
We use this fact to deduce the result for other values of $N$.

For any two choices of $N$, say $N^\ra$ and $N^\rb$, we get different
spaces $V$, $U$, etc.  We distinguish these notationally by using
$V^\ra$ (resp. $V^\rb$) to denote the space $V$ corresponding to 
$N = N^\ra$ (resp. $N^\rb$), and 
likewise for any quantity depending on $N$.
Note that $d = \dim U_1$ is independent of $N$.  

For any fixed $N^\ra$ and $N^\rb$, let
$\phi: V^\ra \to V^\rb$ be the linear map given by
$$
\phi(e_{ij}) =
\begin{cases}
e_{ij'} &\text{where $j'=j+N^\rb-N^\ra$, $S_{ij'} \in R^\rb$}\\
0 &\text{if $S_{ij'} \notin R^\rb$}.
\end{cases}
$$
Let $A$ be the dense open subset of $Gr_d(V^\ra)$,
$$A = \big\{X \in Gr_d(V^\ra)
\ \big|\ X \cap \ker \phi = \{0\} \big\}.$$
Then $\phi$
induces a map $\phi_*:A \to Gr_d(V^\rb)$
$$\phi_*(X) = \mathrm{Image}\ \phi|_X.$$

Observe that $\phi_*(\Udest^\ra) = \Udest^\rb$.
The idea is essentially to show that $A \cap B^\ra \cdot U_1^\ra$ 
is dense in
$\overline{B^\ra \cdot U_1^\ra}$, and that 
$\phi_*(A \cap B^\ra \cdot U_1^\ra) \subset \overline{B^\rb \cdot U_1^\rb}$.
This implies that if $\Udest^\ra \in 
\overline{B^\ra \cdot U_1^a}$, then $\Udest^\rb = \phi_*(\Udest^\ra) \in
\overline{B^\rb \cdot U_1^b}$.  Hence if the result is true for
$N = N^\ra$, then the result will be true for $N=N^\rb$.

We will only prove this in the case where 
$N^\ra=0$ and $N^\rb$ is arbitrary, and in the case where
$N^\ra$ is arbitrary and $N^\rb=0$.  This is enough to give the
result for all $N$.

The case where $N^\ra = 0$ is the easier of the two. 
The map
$$
\left[\begin{matrix}
 A^{n-l \times n-l} & B^{n-l \times l} \\
 0 & C^{l \times l}
\end{matrix} \right]
\in B^\ra 
\mapsto
\left[\begin{matrix}
 A^{n-l \times n-l} & 0 & B^{n-l \times l} \\
 0 & I^{N^\rb \times N^\rb} & 0 \\
 0 & 0 & C^{l \times l}
\end{matrix} \right]
 \in B^\rb 
$$
allows us to view $B^\ra$ as a subgroup of $B^\rb$,
and $\phi_*$ is a $B^\ra$-equivariant inclusion.
Moreover $\phi_*(U_1^\ra) \in \overline{B^\rb \cdot U_1^\rb}$. 
Thus $\phi_*$ takes
$B^\ra \cdot U_1^\ra$ into $\overline{B^\rb \cdot U_1^\rb}$.

%Moreover if 
%${U'}^{(N=N_1)} = \big ({U'_1}^{(N=N_1)}, U_2^{(N=N_1)}, \{0\} \big )$ 
%is transverse then so is the point
%$$\big ( \phi_*({U'_1}^{(N=N_1)}), U_2^{(N=N_2)}, \{0\} \big ) 
%\in Gr(V^{(N=N_2)})^3.$$  
%If the result holds
%for $N=N_1$ {\em and} this latter point
%lies in the $B^3$-orbit closure through $U^{(N=N_2)}$, then the result
%will be true for $N=N_2$ as well.
%
%It suffices to show the result for $N=0$.  
%This is because for $\phi_{0,N_2}$ the map $\phi_*$ is an inclusion 
%of the $B$-orbit closure
%through $U_1^{(N=0)}$ into the $B$-orbit closure though
%$U_1^{(N=N_2)}$.
%
%The result holds for $N$ sufficiently large.  Thus there
%exists a transverse $U'^{(N=N_1)}$ for some $N_1$.
%We now consider the map $\phi_{N_1,0}$.  Note
%that ${U'_1}^{(N=N_1)}$ is in the domain of the induced map $\phi_*$.  
%We show that $\phi_*$ takes a dense subset of $B \cdot U_1^{(N=N_1)}$
%to a dense subset of $B \cdot U_1^{(N=0)}$.  

For the case where $N^\rb = 0$, we consider
the $B$-orbit not through $U_1 \in Gr_d(V)$, but through a lifted 
point $\tilde U_1 \in Gr_d(W)$.  $\tilde U_1$ is defined
in the same way as $U_1$: 
$\tilde U_1 = 
\text{span}\{\tilde e_{ij}\ |\ \text{$S_{ij}$ contains a $1$-token}\}$.
Let $\tilde \phi_*$ be defined analogously to $\phi_*$, taking 
a dense subset of $Gr_d(W^\ra)$ to
$Gr_d(W^\rb)$.
It suffices to show that $\tilde \phi_*$
takes a dense subset of $B^\ra \cdot \tilde U_1^\ra$
to a subset of $B^\rb \cdot \tilde U_1^\rb$.  

Let $L \cong GL(n-l) \times GL(l+N)$ be the 
subgroup of $GL(n{+}N)$ of block diagonal matrices of type $(n-l,l+N)$.
Now $L$ also acts on $W$, and $\tilde U_1$ is fixed by 
$B_- \cap L$.
Since $(B \cap L) \cdot (B_- \cap L)$ is dense in $L$, it follows
that the orbit $B \cdot \tilde U_1 = (B \cap L) \cdot \tilde U_1$ 
is dense in $L \cdot \tilde U_1$.

Thus in fact it suffices to show that $\phi_*$ takes a dense subset of
$L^\ra \cdot \tilde U_1^\ra$ to a subset of 
$L^\rb \cdot \tilde U_1^\rb$.  But this is true, as
\begin{displaymath}
\phi_* \bigg (
 \left[ \begin{matrix}
 A^{n-l \times n-l} & 0 & 0 \\
 0 & B^{l \times N^\ra} & C^{l \times l} \\
 0 & D^{N^\ra \times N^\ra} & E^{N^\ra \times l}
\end{matrix} \right] \cdot
\tilde U_1^\ra \bigg ) = 
\left[ \begin{matrix}
 A^{n-l \times n-l} & 0 \\
 0 & C^{l \times l}
\end{matrix} \right] \cdot \tilde U_1^\rb
\end{displaymath}
whenever both matrices are invertible.

\end{proof}

Although Theorem \ref{thm:geomconverse} is a geometric statement,
our proof ultimately relies on the combinatorics of the 
Littlewood-Richardson rule.  The key non-geometric fact we use
is that $\int_{Gr_l(n)} [Y_\sigma] [Y_\mu] [Y_\nu] \neq 0$ if
and only if there exists a picture between $\lambda_1$ and $\skewtabl$.
It would be an interesting project to find a purely geometric proof 
of this theorem.  The hope would be that a geometric proof of 
Theorem \ref{thm:geomconverse} might allow us to see how to generalise 
some of the results in this paper beyond the Grassmannian.

\subsection*{Acknowledgements}
The author is grateful to Allen Knutson and Stephanie van Willigenburg
for providing feedback and corrections on this paper.


\begin{thebibliography}{KTW}

\bibitem[B]{B} P. Belkale, 
\emph{Geometric proofs of {H}orn and saturation conjectures},
J. Algebraic Geometry \textbf{15} (2006), no. 1, 133-173.

\bibitem[C]{C} I. Coskun, 
\emph{A Littlewood-Richardson rule for two-step flag manifolds},
preprint: {\tt http://www-math.mit.edu/\~{}coskun/}.  

\bibitem[F1]{F1} W. Fulton,
\emph{Young tableaux with applications to representation theory and geometry},
Cambridge U.P., New York, 1997.

\bibitem[F2]{F2} W. Fulton,
\emph{Eigenvalues, invariant factors, highest weights, 
and {S}chubert calculus}, 
Bull. Amer. Math. Soc. \textbf{37} (2000), 209--250.  

\bibitem[FG]{FG} S. Fomin and C. Greene,
\emph{A Littlewood-Richardson miscellany},
European J. Combin. {\bf 14} (1993), no. 3. 191--212.

\bibitem[KTW]{KTW}
A.~Knutson, T.~Tao, and C.~Woodward, 
\emph{The honeycomb model of {${\rm GL}\sb n(\CC)$} tensor products {II}. 
{P}uzzles determine facets of the Littlewood-Richardson cone}, 
J. Amer. Math. Soc. \textbf{17} (2004), no.~1, 19--48 (electronic).

\bibitem[LR]{LR}
D.E. Littlewood and A.R. Richardson, 
\emph{Group characters and algebra},
Philos. Trans. Roy. Soc. London. {\bf 233} (1934), 99--141.

\bibitem[P]{P} K. Purbhoo,
\emph{Vanishing and non-vanishing criteria in Schubert calculus},
Int. Math. Res. Not. (2006), 24590, 1--38.

\bibitem[RW]{RW} J.B. Remmel and R. Whitney,
\emph{Multiplying Schur functions},
J. of Algorithms {\bf 5} (1984), 471--487.

\bibitem[V]{V} R. Vakil,
\emph{A geometric {L}ittlewood-{R}ichardson rule},
to appear in Ann. Math,
arXiv preprint: {\tt math.AG/0302294}.

\bibitem[Z]{Z} A.V. Zelevinsky,
\emph{A generalization of the Littlewood-Richardson rule and the
Robinson-Schensted-Knuth correspondence},
J. Algebra {\bf 69} (1981), no. 1, 82--94.
\end{thebibliography}
\end{document}